\documentclass{amsart}

\usepackage{epsfig}
\usepackage{latexsym}
\usepackage[all]{xy}
\usepackage{amssymb}
\usepackage{amsthm}
\usepackage{amsmath}
\usepackage{amsxtra}

\theoremstyle{plain}

\theoremstyle{remark}

\DeclareMathOperator{\sgnRV}{sgn}

\newcommand{\accaRV}{\`a }

\newcommand{\QRV}{{\mathbb Q}}

\newcommand{\CRV}{{\mathbb C}}

\newcommand{\RRV}{{\mathbb R}}
\newcommand{\ZRV}{{\mathbb Z}}
\newcommand{\NRV}{{\mathbb N}}
\newcommand{\BbRV}{{\mathbb B}}
\newcommand{\TbRV}{{\mathbb T}}

\newcommand{\HcalRV}{{\mathcal{H}}}
\newcommand{\OcalRV}{{\mathcal{O}}}

\begin{document}

\title{Computing central values of $L$-functions}


\maketitle

\section{{}}

 How fast can we compute the value of an $L$-function at the
center of the critical strip?

We will divide this question into two separate questions while
also making it more precise. Fix an elliptic curve $E$ defined
over $\QRV$ and let $L(E,s)$ be its $L$-series. For each
fundamental discriminant $D$ let $L(E,D,s)$ be the $L$-series of
the twist $E_D$ of $E$ by the corresponding quadratic character;
note that $L(E,1,s)=L(E,s)$.

\bigskip
{\bf A.}  \,How fast can we compute the central value $L(E,1)$?

{\bf B.}\, How fast can we compute $L(E,D,1)$ for $D$ in some
interval say
   $a\leq D \leq b$?
\bigskip

  These questions are obviously related but, as we will argue below,
  are not identical.

  We should perhaps clarify what {\it to compute} means. First of all,
  we know, thanks to the work of Wiles and others, that $L(E,s)=L(f,s)$
  for some modular form $f$ of weight $2$; hence, $L(E,s)$, first
  defined on the half-plane $\Re(s)>3/2$, extends to an analytic
  function on the whole $s$-plane which satisfies a functional equation
  as $s$ goes to $2-s$. In particular, it makes sense to talk about the
  value $L(E,1)$ of our $L$-function at the center of symmetry $s=1$.
  The same reasoning applies to $L(E,D,s)$.

  As a first approximation to our question we may simply want to know
  the real number $L(E,D,1)$ to some precision given in advance; but we
  can expect something better. The Birch--Swinnerton-Dyer conjectures
  predict a formula of type
\begin{equation}
\label{waldspurger-fmla} L(E,D,1) = \kappa_D\, m_D^2,
\end{equation}
for some integer $m_D$ and $\kappa_D$ an explicit easily computable
positive constant.  (Up to the usual fudge factors the conjectures
predict that $m_D^2$, if non-zero, should be the order of the
Tate--Shafarevich group of $E_D$.) To compute $L(E,D,1)$ would then
mean to calculate $m_D$ {\it exactly}.

In fact, formulas \accaRV  la Waldspurger have the form
\eqref{waldspurger-fmla} with $m_D$ the $|D|$-th coefficient of a
modular form $g$ of weight $3/2$ which is in Shimura
correspondance with $f$.  The main point of this note is to
discuss informally how explicit versions of such formulas can be
used for problem {\bf B} above.

  Let us also note the interesting fact that $m_D$, being related to
  the coefficient of a modular form, typically does not have a constant
  sign. The significance of the extra information provided by $\sgnRV(m_D)$
  remains a tantalizing mystery.

\section{{}}

There is a standard analytic method to compute $L(E,1)$, which we
now recall. If $E$ has conductor $N$ then the associated modular
form $f$ has level $N$ and
$$
f|_{w_N} = -\varepsilon f,
$$
where $w_N$ is the Fricke involution and $\varepsilon$ is the sign
of the functional equation for $L(E,s)$. Concretely, we have
$$
f\left(\frac i{\sqrt Nt}\right)= \varepsilon
t^2f\left(\frac{it}{\sqrt
     N}\right), \qquad t \in \RRV.
$$
It follows that
$$
\left(\frac{2\pi}{\sqrt{N}}\right)^{-s}\,\Gamma(s) L(f,s) =
\int_0^\infty f\left(\frac{it}{\sqrt     N}\right) t^s \,\frac{dt}
t.
$$
Now break the integral as $\int_0^1 + \int_1^\infty$, make the
substitution $t\mapsto 1/t$ in the first  and use the functional
equation to obtain
$$
\left(\frac{2\pi}{\sqrt{N}}\right)^{-s}\, \Gamma(s) L(f,s) =
\int_1^\infty f\left(\frac{it}{\sqrt
     N}\right) t^s \,\frac{dt} t
+ \varepsilon \int_1^\infty f\left(\frac{it}{\sqrt
     N}\right) t^{2-s} \,\frac{dt} t.
$$
(This is the classical argument to prove the functional equation
of $L(E,s)$ and goes back to Riemann who used it for his zeta
function.)

Now plug in $s=1$ to get
\begin{equation}
   \label{std-fmla}
   \frac{\sqrt{N}}{2\pi}\, L(E,1) = (1+\varepsilon) \sum_{n\geq
   1} \frac{a_n} n \,e^{-2\pi n/\sqrt{N}}
\end{equation}
where $f$ has Fourier expansion
$$
f=\sum_{n\geq 1} a_n \, q^n, \qquad q=e^{2\pi i z}.
$$
Assume for simplicity that $\gcd(N,D)=1$. Then the conductor of
$E_D$ is $ND^2$ and \eqref{std-fmla} applied to $E_D$ yields, more
generally,
\begin{equation}
   \label{std-fmla-D}
   \frac{|D| \sqrt{N}}{2\pi}\, L(E,D,1) = (1+\varepsilon_D) \sum_{n\geq
   1} \left(\frac Dn\right) \frac{a_n} n \,e^{-2\pi n/|D|\sqrt{N}}, 
\end{equation}
with $\varepsilon_D$ the sign in the functional equation of
$L(E,D,s)$.

We know that $|a_n|$ grows no more than polynomially with $n$ (a
straightforward argument gives $|a_n| =O(n)$). It follows that for a
fixed $E$ and varying $D$ we will need to take, very roughly, of the
order of $O(|D|)$ terms in the sum to obtain a decent approximation to
$L(E,D,1)$.  Assuming the Birch--Swinnerton-Dyer conjectures we may
use \eqref{std-fmla-D} to compute $m_D^2$ in \eqref{waldspurger-fmla}
exactly. However, if we know that $m_D$ is the $D$-th coefficient of
some specific modular form (i.e.  we have a formula \accaRV la
Waldspurger) we would get $|m_D|$ but would not be able to recover
$\sgnRV(m_D)$.

Using this method to compute, say, $L(E,D,1)$ for $|D|\leq X$
would take time of the order of $O(X^2)$.  We will see below that
using formulas of type \eqref{waldspurger-fmla} we can reduce this
to $O(X^{3/2})$ for at least some fraction of such $D$'s.

\section{{}}
Before tackling $L(E,D,1)$ let us consider the case of the special
value of an Eisenstein series of weight $2$ (as opposed to a cusp
form as we have for $L(E,D,1)$). What follows is meant only as an
illustration of the general case.

Let the $L$-function be $L(\left(\frac D
   \cdot\right),s-1)L(\left(\frac D \cdot\right),s)$ with $D<0$ the
discriminant of an imaginary quadratic field $K$.  Its value at
$s=1$ is essentially $h(D)^2$, where $h(D)$ is the class number of
$K$, and we find an analogue of \eqref{waldspurger-fmla} with
$h(D)$ playing the role of $m_D$.  There are many excellent
algorithms for computing the class number $h(D)$ (see for example
\cite{C} chap.~$5$). Unfortunately, these do not obviously
generalize to the calculation of $m_D$. The main reason for this
is that the class group of $K$ is easy to describe (both its
elements and the group operation) in terms of binary quadratic
forms, whereas its elliptic analogue, the Tate-Shafarevich group
of $E_D$, is notoriously intractable.

The standard analytic method of the previous section yields the
following formula (which was known to Lerch, see \cite{D}
vol.~III, p.~171)
\begin{equation}
\label{lerch-fmla} h(D)^2=\frac{w_D^2\sqrt{|D|}}{2\pi} \sum_{n\geq
1}\left(\frac D
   n\right) \frac{\sigma(n)} n \,e^{-2\pi n/|D|},
\end{equation}
where $w_D$ is the number of units in $K$ and
$\sigma(n):=\sum_{d\mid
   n}d$ is the divisor sum function. Again, we need to take,
roughly, $O(|D|)$ number of terms in the sum to obtain a
reasonable approximation of the left hand side. In this case, we
in fact have an {\it exact} formula requiring $D$ terms, namely,
Dirichlet's class number formula
\begin{equation}
\label{dirichlet-fmla} h(D)=-\frac {w_D}{2|D|}
\sum_{n=1}^{|D|-1}n\left(\frac D
   n\right).
\end{equation}
Neither one of these formulas is, however, particularly useful for
computing $h(D)$ in practice. On the other hand, it may be worth
pointing out that  similar arguments yield the formula \cite{D}
vol.~III, p.~153.
$$
h(D)= w_D \sum_{n\geq 1} \left(\frac D
   n\right)\frac 1{1-(-1)^ne^{\pi n/\sqrt{|D|}}}, \qquad D\equiv 5
\bmod 8,
$$
with the number of necessary steps now reduced to the order of
$O(\sqrt{|D|})$. (Analogous formulas can be given for $D$ in other
congruence classes modulo $8$.)

To make the connection to the general case of computing $L(E,D,1)$
that we are considering we mention two other possible approaches
to computing $h(D)$ that do generalize.

\medskip
\noindent {\bf (I)} The first is to follow Gauss and realize ideal
classes of $K$ as classes of primitive, positive definite binary
quadratic forms of discriminant $D$. Each class has a unique
representative $Q=(a,b,c)$ in the standard fundamental domain
(what is known as a {\it reduced
   form}) and we can simply enumerate these. A straightforward
algorithm is as follows: run over values of $b$ with $b\equiv
D\bmod 2$ and $0\leq b \leq \sqrt{|D|/3}$; for each $b$ decompose
$(b^2-D)/4$ as $ac$ with $0<a\leq c$. Add one or two to the total
count as the case may be if $\gcd(a,b,c)=1$.

Though this algorithm also takes time $O(|D|)$ the constant of
proportionality is very small making the algorithm quite
practical. An important point to notice for our purpose, however,
is that if we wanted to compute $h(D)$ for $0\leq |D|\leq X$ we
may simply run over all triples $a,b,c$ of size at most
$\sqrt{X/3}$ checking the necessary conditions on $(a,b,c)$ for it
to be a reduced form. In this way we obtain an algorithm which
will run in time $O(X^{3/2})$.

\medskip
\noindent {\bf (II)} The second approach is again to follow Gauss
but in a different direction. He proved that $h(D)$ is related to
the number of representations of $|D|$ as a sum of three squares.
One  precise form of this relation is the following identity (see
\cite{G} p.177)
\begin{equation}
\label{gauss} \tfrac12 \sum_{x\equiv y\equiv z \bmod 2}
q^{x^2+y^2+z^2}= \tfrac12+ 12\sum_{D}H_2(D)\,q^{|D|}
\end{equation}
where $D$ runs through all negative discriminants (i.e. $D<0$ and
$D\equiv 0,1 \bmod 4$), and $H_2$ is a variant of the Hurwitz
class number (see \cite{G}, page 120). (For us it suffices to know
that it is related to $h(D)$; for example for $D \equiv 5 \bmod 8$
a fundamental discriminant we have $H_2(D)=h(D)$.)

There are sophisticated techniques for computing the coefficients
of the left hand side, such as {\it convolution} which uses the
fast Fourier transform to compute products of $q$-series. But even
a simple enumeration of the lattice points $x^2+y^2+z^2\leq X,
x\equiv y\equiv z \bmod 2$ would again take time $O(X^{3/2})$.

The two approaches (I) and (II) are of course related; they amount
to {\it counting} (in an appropriate sense) the number of
representations of $D$ by a certain ternary quadratic form. In
case (I) we count the number of solutions to $b^2-4ac=D$ up to
$SL_2(\ZRV)$-equivalence; in (II), the number of solutions to $|D|
= x^2+y^2+z^2$ with $x\equiv y\equiv z \bmod 2$. Note the crucial
difference that the ternary quadratic form involved is indefinite
in case (I) and positive definite in case (II).

A more geometrical point of view is to think that we are dealing
with {\it Heegner points}. In case (I) we may associate to a
primitive positive definite binary quadratic form $Q=(a,b,c)$ the
point $z_Q=(-b+\sqrt{D})/2a$ in the upper half plane $\HcalRV$.
The respective actions of $SL_2(\ZRV)$ on forms and $\HcalRV$ are
compatible; hence, the class of $Q$ determines a unique (Heegner)
point in $SL_2(\ZRV)\backslash \HcalRV$ of discriminant $D$.

  It is a bit less intuitive how to think of Heegner points in case
  (II) but this was worked out by Gross \cite{G}. The main ingredient is a
  positive definite quaternion algebra $B$ over $\QRV$ ramified, say, at
  $\infty$ and a prime $N$. Pick a maximal order $R$ of $B$ and let
  $I_1,\ldots,I_n$ be representatives for the (left) ideal classes of
  $R$. Let $R_i$ be the right order of $I_i$ for $i=1,\cdots,n$.

  Fix an imaginary quadratic field $K$ of discriminant $D$. Then we can
  think of a Heegner point of discriminant $D$ (what Gross calls a {\it
    special point}) as an (optimal) embedding of the ring of integers
  $\OcalRV_K$ into some $R_i$.  Eichler has proved that the total number of
  such points, each counted up to conjugation by $R_i^\star$, is
  $(1-\left(\frac D N\right))h(D)$. (In fact, the situation is quite
  analogous to that of case (I) if we take the {\it indefinite} algebra
  $B=M_2(\QRV)$ and $R=M_2(\ZRV)$.)

  For example, if $N=2$ then the algebra $B$ is the usual Hamilton
  quaternions and we may pick $R$ to be the order discovered by Hurwitz
  (in standard notation)
$$
R=\ZRV+\ZRV i+\ZRV j+\ZRV\tfrac12(1+i+j+k).
$$
In this case there is only one class of left $R$-ideals
represented by $R$ itself. Hence a Heegner point is an embedding
$\phi: \OcalRV_K \rightarrow R$.

How do we find such embeddings?  The main thing we need is a $w\in
R$ with $w^2=D$. Such a quaternion, because $D$ is a scalar,
necessarily has trace $t(w)= 0$ and norm $n(w)= -D$ and
conversely. Elements of trace $0$ in $R$ form a rank $3$ lattice
and hence $n(w)=-D$ is a representation of $D$ by a certain
ternary quadratic form associated to $R$. A few congruence
conditions are needed to actually produce an optimal embedding
$\phi$ out of $w$ but the upshot is that the problem becomes one
about representations of $-D$ by ternary quadratic forms. For
example, in the case $N=2$ Eichler's count of embeddings can be
completely encoded in the identity \eqref{gauss}; the presence of
the factor $12$ in that formula is due to the fact that this is
the order of $R^\star/\pm1$. More details on this setup are given
below in \S 4 (II).

\section{{}}
We now return to the main case of computing $L(E,D,1)$ and
describe analogues of cases (I) and (II) of the previous section.
These analogues are the remarkable results of Gross and Zagier.

\medskip
\noindent {\bf (I)} Let us assume for simplicity that $E$ has
conductor a prime $N$, sign of the functional equation equal to
$-1$, $L'(E,1)\neq 0$, and $E(\QRV)=\langle P_0\rangle$. If $f$ is
the weight $2$ eigenform associated to $E$ then we get a map
\begin{equation}
   \label{modular-map}
\begin{array} {cccc}
   \Phi:  \quad & X_0(N) & \longrightarrow & \CRV/L\\
                & z      & \mapsto         & 2\pi i
                \int_{i\infty}^z f(u)\,du
\end{array}
\end{equation}
where $X_0(N)$ is the modular curve of level $N$ and $L\subset
\CRV$ is a certain lattice of periods of $f$. It is known that
$\CRV/L=E'(\CRV)$ for some elliptic curve  $E'/\QRV$ isogenous to
$E$. Since the $L$-function is unchanged by isogenies we may
assume without loss of generality that $E'=E$.

Let $K$ be an imaginary quadratic field of discriminant $D<-4$ in
which $N$ splits. Choose $b_*\in \ZRV$ such that $b_*^2\equiv D
\bmod N$; this is possible by the assumption that $N$ splits in
$K$.  Note also that $N$ does not divide $b_*$. We want to
consider Heegner points on $X_0(N)$ of discriminant $D$. To define
them concretely choose representatives $Q=(a,b,c)$ of the $h(D)$
classes of binary quadratic forms with $N\mid a$ and $b\equiv b_*
\bmod N$. (For example, start with representatives $(a,b,c)$ with
$\gcd(a,N)=1$ and compose them with the fixed form
$(N,b_*,(b_*^2-D)/2N)$.)

Then $z_Q:=(-b+\sqrt{D})/2/a)\in X_0(N)$ is well defined and
$P_D:=\sum_Q \Phi(z_Q)\in E(K)$. Moreover, complex conjugation
fixes $P_D$, by the assumption on the sign of the functional
equation. Hence $P_D$ actually is in $E(\QRV)$ (and is independent
of the choice of $b_*$).

One consequence of the results of Gross--Zagier is the following
\cite{Z},\cite{GZ}, \cite{GKZ}. By our assumption on $E(\QRV)$ we
have $P_D=m_DP_0$ for some $m_D\in \ZRV$ and hence
\begin{equation}
   \label{gross-zagier}
L(E,D,1) = \kappa_D m_D^2;
\end{equation}
where $\kappa_D$ is an explicit easily computable positive constant;
i.e.  we have a formula of type \eqref{waldspurger-fmla}.

Usually one regards the Gross-Zagier formula as a way to compute a
rational point $P_D$ on $E$ whose height is given in terms of
$L(E,D,1)L'(E,1)$ and hence obtaining, when this value does not
vanish, a confirmation of the predictions of the
Birch--Swinnerton-Dyer conjecture. Here, instead, we are taking
the point of view that the points of $E(\QRV)$ are known and use
the Gross--Zagier formula as a means to computing $L(E,D,1)$.

To calculate $m_D$ in practice it is better to work on the 
$E(\CRV)=\CRV/L$ model of $E$ rather than, say, a Weierstrass
equation. Let $z_0\in \CRV$ represent, modulo $L$, the point $P_0
\in E(\QRV)$. We first compute an approximation to
$$
z_D:=\sum_Q\sum_{n\geq 1}  \frac{a_n} n \, e^{2\pi i n z_Q}.
$$
Then we solve the linear equation below for integers $n_1$ and
$n_2$
$$
z_Q = m_Dz_0 +n_1\omega_1+n_2\omega_2,
$$
where $\omega_1,\omega_2$ are a basis for $L$. (In fact,
multiplying by $2$ if necessary, we may assume that $\omega_1\in
\RRV$ and $\omega_1\in i\RRV$ and hence by taking real parts solve
only a three term equation instead.)

The result is a practical and reasonably efficient algorithm for
computing $m_D$. The number $m_D$ is the $D$-th Fourier
coefficient of a weight $3/2$ modular form $g$ of level $4N$ which
is in Shimura correspondence with $f$.  It is interesting that we
can compute the Fourier coefficients of $g$ directly without any
knowledge of the whole vector space of modular forms in which $g$
lies; though we do, of course, start by knowing $f$ itself. (We
have only described the calculation for certain $D$'s but there is
analogous way to get all coefficients.)

Together with my student Ariel Pacetti we implemented the above
algorithm  in GP. The corresponding routines can be found at

\begin{verbatim}
http://www.ma.utexas.edu/users/villegas/cnt/
\end{verbatim}

\noindent under Heegner points.

Here is a sample example. Let $E$ be the curve $y^2+y=x^3-x$ of
conductor $N=37$ (this is the elliptic curve over $\QRV$ of
positive rank with smallest conductor). This case was described in
detail in \cite{Z}. It is known that
$E(\QRV)=\langle(0,0)\rangle$.

\begin{verbatim}
? e=ellinit([0,0,1,-1,0]);  anvec=ellan(e,5000);
? for(d=5,100, if(isfundamental(-d) && kronecker(-d,37)==1,
          print(-d," ",ellheegnermult(e,-d,[0,0],0,anvec)[1])))

-7  -11 -40 -47 -67  -71 -83  -84 -95
 1   -1  -2   1  -6   -1   1    1   0 

\end{verbatim}
The first row is $D$, the second $m_D$ (for typographical reasons we
transposed the actual GP output). These values agree, fortunately,
with Zagier's \cite{Z} formula (28) up to a global negative sign.

In our implementation at least the algorithm is not that well
suited for computing $L(E,D,1)$ for all $D<0$ and $|D| <X$ for
very large $X$; for this, it would be better to adapt (see \S 5)
the ideas of (II) below but these have not been fully implemented
as yet.

\medskip
\noindent {\bf (II)}   Let $B$ over $\QRV$ be the (unique up to
isomorphism) positive definite quaternion algebra ramified at
$\infty$ and a prime $N$. Pick a maximal order $R$ of $B$ and let
$I_1,\ldots,I_n$ be representatives for the (left) ideal classes
of $R$. Let $R_i$ be the right order of $I_i$ for $i=1,\cdots,n$.
The class number $n$ of $R$, in contrast with $h(D)$, has a simple
formula and is roughly of size $N/12$.

For example, if $N\equiv 3 \bmod 4$ we can describe $B$ as the
algebra over $\QRV$ with generators $i,j$ such that
$i^2=-1,j^2=-N$ and $ij=-ji$. Also in this case we can take
$R=\ZRV+\ZRV i+\ZRV\tfrac12(1+j)+\ZRV i\tfrac12(1+j)$.

There are various ways to compute representatives $I_1,\cdots,
I_n$ of the ideal classes (for algorithms for quaternion algebras
see \cite{P}). If $N\equiv 3 \bmod 4$ there is an algorithm which
is completely analogous to that of Gauss \S3 (I) for binary
quadratic forms. It exploits the fact that our choice of $R$ has
an embedding of $\ZRV[i]$ and hence allows us to view $R$-left
ideals as rank $2$ modules over $\ZRV[i]$; then classes of
$R$-ideals correspond to classes of positive definite binary
Hermitian forms over $\ZRV[i]$ of discriminant $-N$. Instead of
$\HcalRV$ we now need to work on hyperbolic $3$-space where, as it
turns out, the action of $SL_2(\ZRV[i])$ has a very simple
fundamental domain. This yields an algorithm which is almost
verbatim that of Gauss for binary forms over $\ZRV$. Details can
be found in \cite{RV}.

  For example, if $N=11$ then there are two classes of positive
  definite binary Hermitian forms of discriminant $-11$ over $\ZRV[i]$;
  namely, $(1, 1, 3)$ and $(2, 1 + 2i, 2)$ corresponding to the two
  ideals
$$
I_0:=R=\ZRV+\ZRV i+\ZRV\tfrac12(1+j)+\ZRV i\tfrac12(1+j)
$$
and
$$
I_1:=2\ZRV+\ZRV 2i+\ZRV\tfrac12(1+2i+j)+\ZRV i\tfrac12(1+2i+j)
$$
representing the $n=2$ classes of left $R$-ideals.

Let $V_\QRV$ be the $\QRV$ vector space of functions on the set
$\{I_0, \ldots, I_n\}$. For each $m\in \ZRV_{\geq 0}$ there is an
operator $B(m)$ acting on $V_\QRV$, the Brandt matrix of order
$m$, which encodes the number of representations of $m$ by certain
quaternary quadratic forms (see \cite{G} (1.4)). Let $\BbRV$ be
the algebra generated over $\ZRV$ by all the $B(m)$; it is
commutative and $\BbRV\otimes_\ZRV \QRV$ is semisimple.

On the other hand, we have the space $M_\CRV$ of modular form of
weight $2$ on $\Gamma_0(N)$ (known to be of dimension $n$) and the
Hecke operators $T_m$ acting on $M_\CRV$. Let $\TbRV$ be the
algebra spanned by the $T_m$ over $\ZRV$ ; like $\BbRV$ it is
commutative and $\TbRV\otimes_\ZRV \QRV$ is semisimple.  This
algebra preserves the $\QRV$ vector space $M_\QRV\subset M_\CRV$
of dimension $n$ consisting of those modular forms in $M_\CRV$
with Fourier coefficients in $\QRV$.

These two setups are closely related and indeed we have a special
case of the Jacquet--Langlands correspondence. Eichler proved that
$T_m$ and $B(m)$ have the same trace for all $m\in \NRV$. Hence,
by semisimplicity of the algebras the map $T_m\mapsto B(m)$
induces a ring isomorphism $\TbRV \simeq \BbRV$. It follows that
eigenspaces of $V_\QRV$ and $M_\QRV$, under the action of $\BbRV$
and $\TbRV$ respectively, correspond to each other. Since we also
have multiplicity one these eigenspaces are one-dimensional.

In conclusion, given $f=\sum_{n\geq 0} a_n q^n \in M_\CRV$ an
eigenform for all Hecke operators $T_m$ (so that $T_m f = a_m f$)
there is an $e_f \in V_\QRV\otimes_\ZRV K$ unique up to scalars
such that $B(m)e_f=a_me_f$. (Here $K$ denotes the field
$\QRV(a_0,a_1,\ldots)$ generated by the Fourier coefficients of
$f$.)

In fact, this correspondence gives an efficient way to compute
Fourier coefficients of eigenforms in $M_\CRV$ (see \cite{P}).  An
implementation of the corresponding algorithms can be found in the
above mentioned website (under {\tt qalgmodforms}). Here is a
sample GP session.

\begin{verbatim}

? R=qsetprime(11);

? brandt(R,2)~

[1 3]

[2 0]

? brandt(R,3)~

[2 3]

[2 1]

\end{verbatim}
The first line defines $R$ as a maximal order in the algebra
ramified at $11$ and $\infty$; the others compute the
corresponding Brandt matrices. We find that these matrices have
two eigenvectors: $e_E=(1/2,1/3)$ and $e_f=(-1,1)$ corresponding to an
Eisenstein series and a cusp form, respectively.

The above implementation is intended for small to medium scale
computations. For large scale computations one should use the {\it
  graph method} ideas of Mestre and Oesterl\'e \cite{M}, which exploit
the sparse nature of the Brandt matrices.

Now following Gross we show how to associate a modular form of
weight $3/2$ to an eigenvector $e_f$. Let $R_i$ be the right order
of $I_i$ and let $L_i\subset R_i$ be the ternary lattice defined
by
$$
L_i: \qquad w \in R_i, \qquad t(w)=0, \qquad w \in \ZRV \bmod
2R_i.
$$
Let $g_i$ be the corresponding theta series
$$
g_i(\tau) := \tfrac12 \sum_{w\in L_i}q^{n(w)}, \qquad q=e^{2\pi
i\tau}.
$$
Gross \cite{G} prop. 12.9 describes precisely how the $D$-th
coefficient $a_i(D)$ of $g_i$ relates to the optimal embeddings of
imaginary quadratic orders of $Q(\sqrt{D})$ into $R_i$.

These theta series are modular forms of weight $3/2$ and level
$4N$ and, in fact, belong to a certain subspace $U$ defined by
Kohnen. This subspace is determined by the condition that the
coefficient of $q^d$ of a form should be zero unless $D:=-d$ is a
discriminant, i.e., $D\equiv 0,1 \bmod 4$, and $\left(\frac D
N\right)\neq 1$. The weight $3/2$ Hecke operators $T_{m^2}$
preserve $U$.

Define
$$
g:=\sum_i e_f(i)\,g_i = \sum_D m_D \, q^{|D|} \in U.
$$
This form is identically zero if the sign in the functional
equation of $f$ is $-1$. If $g$ is non-zero it is a modular form
in Shimura correspondence with $f$; i.e., $T_{m^2}g = a_m g$,
where $T_m f = a_m f$. Moreover, we have the  Waldspurger formula
\cite{G} 13.5
\begin{equation}
   \label{gross-waldspurger}
  L(f,1)L(f\otimes \chi_D,1) =\kappa_f\frac{\delta_D}{\sqrt{|D|}}\, m_D^2,
\end{equation}
where $D$ is a fundamental discriminant with $\left(\frac D
   N\right)\neq 1$, $\chi_D$ is the associated quadratic character,
$\kappa_f >0$ is a constant depending only on $f$ and
$\delta_D:=2$ if $N \mid D$ and $\delta_D:=1$ otherwise.

Finally, let $E/\QRV$ be an elliptic curve of prime conductor $N$
and sign $+1$ in its functional equation. Let $f$ and $g$ be the
corresponding modular forms of weight $2$ and $3/2$ respectively
as above. Then if $L(f,1)\neq 0$ we obtain from
\eqref{gross-waldspurger} a formula of type
\eqref{waldspurger-fmla} with $m_D$ the Fourier coefficient of
$g$. As in \S 3 (II) to compute $m_D$ for $|D|<X$ we could run
through all $w\in L_i$ with $n(w)\leq X$ whose total number is
$O(X^{3/2})$. Again various computational techniques could also be
used to speed up the calculation of $m_D$. Note that in any case
all computations are done with integer arithmetic.

Tables of $m_D$'s for several curves and the routines to compute
them can be found at G. Tornar\'{\i}a's website

\begin{verbatim}
http://www.ma.utexas.edu/users/tornaria/cnt/
\end{verbatim}

\noindent among other goodies (an interactive version of Cremona's
tables of elliptic curves and an interactive table of ternary
quadratic forms).

\section{{}}
We conclude with some remarks about the general situation.

\medskip
\noindent {\bf 1. \,} It follows from \eqref{gross-waldspurger}
that if $L(f,1)=0$ then the form $g$ vanishes identically. In this
case we naturally need to do something else.

In \cite{MRVT} we work out an extension of Gross's work
introducing an auxiliary prime $l$; the theta series $g_i$, for
example, are modified by introducing an appropriate weight
function. The complexity of algorithms only increase by a factor
essentially proportional to $l$.

\medskip
\noindent {\bf 2. \,} If the level $N$ is not prime but
square-free the situation is not too different from the one
described above. The downside is that $L(E,D,1)$ can be computed
this way only for a certain fraction of $D$'s (determined by local
conditions). One needs to consider a quaternion algebra $B$
ramified at $\infty$ and at primes $l\mid N$ for which the
Atkin-Lehner involution acts as $f|_{w_l}=-f$ and an Eichler order
in $B$ of level the product of the remaining primes factors of
$N$.

\medskip
\noindent {\bf 3. \,} If the level is not square-free things
become quite a bit more complicated; for example, the algebra
$\BbRV$ of Brandt matrices typically does not act with
multiplicity one and some modular forms are simply missing. The
arithmetic of the corresponding orders, which are no longer
Eichler orders in general, also becomes more involved and,
moreover, one needs to consider two types of orders: one for the
weight $2$ side and another for the weight $3/2$ side; see
\cite{PRV}, \cite{PTa}, \cite{PTb} some work on this case.

\medskip
\noindent {\bf 4. \,} To compute twists $L(f\otimes \chi_l,1)$ by
{\it real} quadratic fields $\QRV(\sqrt{l})$ one may consider a
twist $f_D:=f\otimes \chi_D$ by an auxiliary imaginary quadratic
field $Q(\sqrt{D})$ and find a formula of type
\eqref{waldspurger-fmla} for $L(f_D\otimes \chi_{Dl},1)$. The form
$f_D$ typically does not have square-free level so several
corresponding difficulties ensue, see \cite{PTb}.

\medskip
\noindent {\bf 5. \,} Forms of higher weight can also be handled
using quaternion algebras by introducing harmonic polynomials as
weight functions for the theta functions (both for the ideals
$I_i$ corresponding to forms of weight $2+2r$ and for the ternary
lattices $L_i$ corresponding to forms of weight $3/2+ r$) see
\cite{H}, \cite{RT}.

\begin{flushleft}
\mbox{}\\
\mbox{}\\
Department of Mathematics\\ University of Texas at Austin,\\ TX
78712 USA\\
\mbox{}\\
\mbox{}\\
villegas@math.utexas.edu
\end{flushleft}

\end{document}